\documentclass[draft]{amsart}
\usepackage{enumitem,mathtools,xcolor}
\usepackage{amsthm}
\usepackage{amsmath}
\usepackage[mathscr]{eucal}
\theoremstyle{plain}
\newtheorem{theorem}{Theorem}

\newtheorem{lemma}{Lemma}
\newtheorem{proposition}{Proposition}

\theoremstyle{definition}

\newtheorem{assumption}[theorem]{Assumption}

\theoremstyle{remark}
\newtheorem{remark}{Remark}

\DeclareMathOperator{\diag}{diag}

\allowdisplaybreaks[3]

\begin{document}

\title{Sixth-order Birkhoff regular problems}

\address{
 School of Mathematics,
 University of the Witwatersrand,
 Johannesburg 2050, South Africa
}
\address{
 National Institute for Theoretical and Computational Sciences (NITheCS),
 South Africa
}

\author{Nokukhanya Thandiwe Mzobe}
\email{thandiwe.mzobe1@wits.ac.za}
\author{Bertin Zinsou}
\email{bertin.zinsou@wits.ac.za}

\begin{abstract}
Asymptotics of the eigenvalues can always be derived for self-adjoint boundary value problems. However, they can also be derived for boundary value  problems that fail to be self-adjoint provided that they are Birkhoff regular.
A regular sixth-order differential equation that depends quadratically on the eigenvalue parameter $\lambda$  is
considered with classes of separable boundary conditions independent of $\lambda$  or depending linearly on $\lambda$. Conditions are
given for the problems to be Birkhoff regular.
\end{abstract}

\maketitle

\section{Introduction}

\par The development of spectral theory is one of the most significant chapters in the history of  mathematics. Spectral theory came from an attempt to provide a mathematical framework for understanding various physical phenomena and later developed into a field that has greatly enriched mathematics as a whole. It  has been studied through many qualitative and quantitative techniques such as the Sturm-Liouville theory, separation of variables, Fourier and Laplace transforms, perturbation theory, eigenfunction expansions, variational methods, stochastic analysis and numerical methods including
finite elements to mention a few.

 Spectral theory remains a vital field of on-going investigation in Mathematics.  One of the most notable contributions  of spectral theory is the theory of linear  differential systems and their generalizations. The spectral theory of ordinary differential equations is the part of spectral theory concerned with the determination of the spectrum and eigenfunction expansion associated with a linear ordinary differential equation. Ordinary differential equations have important applications  in the study of many problems in the natural sciences, physics, engineering and aeronautics. These applications can be investigated through linear operators.


Sturm-Liouville operators  operators have provided a constant source of new ideas and problems in the spectral theory of operators  while higher-order differential operators are experiencing slow but steady developments. Higher order linear differential esectionquations occur in applications with or without the eigenvalue parameter in the boundary conditions. For a comprehensive study of problems with dependent boundary conditions, we refer the reader to
\cite{moller2006spectral,moller2011self,moller2012spectral,moletsane2017self,pivovarchik2001inverse,moller2013sixth,moller2015self,moller2017asymptotics,binding1994sturm,aliyev2012spectral,marletta2003pencils}. The mathematical model of these eigenvalue problems leads to differential operators which depend quadratically on the eigenvalue parameter. These problems are recognized as operator polynomials, a synonym for operator pencils, and their operator representation is given by
\begin{equation}\label{1}
L(\lambda)=\lambda^{2}M-i\lambda K-A
\end{equation}
in the Hilbert space $H=L_{2}(I) \ \oplus \mathbb{C}^{l}$, where $ I $ is an interval, $ l $ is the number of boundary conditions depending on the eigenvalue parameter, $M, K$ and $A$ are coefficient operators.

 In the theory of boundary value problems for differential equations, a large part is played by differential operators with regular boundary conditions. Birkhoff was the first to introduce regular boundary conditions for ordinary differential operators at the beginning of the 20th century \cite{birkhoff1908boundary}. One of the most significant consequences of Birkhoff regularity of a boundary eigenvalue problems is  that the determinant of a suitable characteristic matrix function is bounded away from zero for $\lambda$ in the union of circles $\Gamma_{\nu}, ( \nu\in\mathbb{N})$ with centres at zero and radii $\rho_{\nu}$ which tend to infinity as $\nu$ tends to infinity \cite[page 129]{mennicken2003non}.

 In \cite{birkhoff1908boundary}, Birkhoff proved a theorem on eigenfunction expansion and an estimate for Green’s function of the class of regular boundary conditions satisfying some determinental inequalities. Stone \cite{stone1926comparison} proved that Birkhoff expansion is in a sense equivalent to the Fourier series. This result was partially justified by Salaff in his remarkable paper \cite{salaff1968regular}, where he proved that self-adjointness implies Birkhoff regularity for even order operators.  Birkhoff regular boundary conditions of a small transversal vibrations of a homogeneous beam compressed or stretched were investigated in \cite{moller2017asymptotics,moller2006spectral,moller2012spectral,zinsou2016fourth,zinsou2021}. Most of the studies conducted on Birkhoff regular boundary conditions  focus on fourth order eigenvalue problems. Furthermore, it appears that there is almost no work on sixth-order Birkhoff regular problems. As a result, we believe it is still an open problem and it is imperative to initiate its investigation. 

We  investigate  sixth-order eigenvalue problems with Birkhoff regular boundary conditions. The application of sixth-order differential operators is found in mathematical models of vibrations of curved arches, see \cite{auciello1994free}. Numerical methods and other techniques for investigating solutions of sixth-order boundary value problems  can be found in \cite{boutayeb1992numerical,twizell1988numerical,mohyud2009variation,loghmani2007numerical}.
 Moller and Zinsou  \cite{moller2013sixth}  investigated the spectral properties of a sixth-order problem described by the operaor pencil 
%
%
\begin{equation}\label{Pencil}
L(\lambda)=\lambda^{2}M-i\lambda K-A,
\end{equation}
where $M$,  $K$ and $A$ are self-adjoint operators. 

We consider the sixth order boundary value problem given by the differential equation
\begin{equation}\label{AB}
-y^{(6)}+(g_{2}y'')''-(g_{1}y')'+g_{0}y=\lambda^{2}y,\\
\end{equation}
where $g_{m}\in C^{m}[a,b]$, for $m=0,1,2$, are real valued, and a class of separable boundary conditions
\begin{equation}\label{3}
B_{j}(\lambda)=0 \qquad j=1,2,3,4,5,6,
\end{equation}
 The boundary conditions \eqref{3} are independent of $\lambda$ or dependent on $\lambda$ linearly. The quadratic operator pencil associated with the problems under consideration is given by \eqref{Pencil},
where the main operator A is not necessarily self-adjoint. We give conditions for the problem \eqref{AB}, \eqref{3} to be Birkhoff regular. In a forthcoming paper we will investigate the asymptotics of the eigenvalues for which the sixth order boundary value problems are Birkhoff regular.

 We introduce the eigenvalue problems under consideration in Section \ref{Sec2}. In Section \ref{Sec3}, we provide the setting needed to conduction our investgation.  Finally, in Section \ref{Sec4}, we provide the results of our investigation.


\section{Sixth order eigenvalue problem} \label{Sec2}
 
We consider  the sixth-order eigenvalue problem
\begin{gather}\label{3.1}
\displaystyle -y^{(6)}+\sum\limits_{r=0}^2\left(g_ry^{(r)}\right)^{(r)}=\lambda^{2}y,\\
\label{3.2} B_{j}(\lambda)y=0, \quad j=1,2,3,4,5,6,
\end{gather}
where $y\in W_{2}^{6}(a,b)$,  $g_{r}\in C^{r}[a,b]$, for $r=0,1,2$, are real valued functions, and \eqref{3.2} are separated boundary conditions lineary dependent or independent on $\lambda$.  We assume that

\begin{equation}\label{3.3}
B_{j}(\lambda)y=\sum_{k=0}^{p_{j}}\alpha_{j,k}y^{(k)}(a_{j})+i\lambda\sum_{k=0}^{q_{j}}\beta_{j,k}y^{(k)}(a_{j}),\\
\end{equation}

\noindent where   $p_{j},q_{j}\in\lbrace-\infty,0,1,2,3,4,5\rbrace$,   at least one of the numbers $p_{j},q_{j}\neq-\infty$, $j\in\lbrace{1,2,3,4,5,6\rbrace}$. We set $a_{j}=a$ for $j=1,2,3$ and $ a_{j}=b$ for $j=4,5,6$. We set
\begin{equation}\label{alphabeta}
\alpha_{j,p_{j}}=1 \text{ if }p_{j}\neq-\infty \quad \text{ and } \quad \beta_{j,q_{j}}\neq0 \text{ if }  q_{j}\neq-\infty.\\
\end{equation}
Define
\begin{gather*}
\Theta_{1}=\lbrace s\in\lbrace1,2,3,4,5,6\rbrace: B_{s}(\lambda) \ \text{depends on } \lambda\rbrace,\ \
\Theta_{0}=\lbrace1,2,3,4,5,6\rbrace\backslash\Theta_{1},\\	\Theta_{1}^{a}=\Theta_{1}\cap\{1,2,3\},\ \ \Theta_{1}^{b}=\Theta_{1}\cap\{4,5,6\},
\end{gather*}
and 
\begin{gather}\label{3.4}
\Lambda=\{s\in \{1,2,3,4,5,6\}:p_{s}>-\infty\}, \ \  \Lambda^{a}=\Lambda\cap\{1,2,3\}, \ \ \Lambda^{b}=\Lambda\cap\{4,5,6\}.
\underline{}\end{gather}

\begin{assumption}\label{3.1.1}
	\textit{ We assume that the numbers $p_{s}$ for $s\in\Lambda^{a}$,  $q_{j}$ for $j\in \Theta_{1}^{a}$ are distinct and that the numbers $p_{s}$ for $s\in\Lambda^{b}$,  $q_{j}$  for $j\in\Theta_1^{b}$ are distinct.}
\end{assumption}

\noindent  The above assumption means for any pair $(r,a_{j})$, the term $y^{(r)}(a_{j})$ appears at most once in the boundary conditions \eqref{3.2}.

Denoting the collection of boundary conditions in  \eqref{3.2}  by $U$, we  define the following operators associated with $U$:
\begin{gather}\label{3.1.1*}
U_{r}y=\left(\sum\limits_{k=0}^{p_{j}}\alpha_{_{j,k}}y^{(k)}\right)_{j\in\Theta_r},\  r=0,1 \text{ and } \ V_{1}y=\left(\sum\limits_{k=0}^{q_{j}}\beta_{j,k}y^{(k)}\right)_{j\in\Theta_{1}}, \,
\end{gather}
 where $y\in W_{2}^{6}(a,b)$ and $W_{2}^{6}(a,b)$ is the Sobolev space of order 6 on the interval $(a,b)$.

Putting $l:=|\Theta_{1}|$, 
 we consider the linear operator $A(U), K \text{ and } M$ in the space $L_{2}(a,b)\oplus\mathbb{C}^{l}$ with domains
\begin{align*}
&\mathscr{D}(A(U))=\left\{\widetilde{y}=\begin{pmatrix}
y\\
V_{1}y
\end{pmatrix}:y\in W_{2}^{6}(a,b), U_{0}y=0\right\},\\
&\mathscr{D}(K)=\mathscr{D}(M)= L_{2}(a,b)\oplus\mathbb{C}^{l},\\ 
\end{align*}
where
\begin{gather*}
(A(U))\widetilde{y} = \begin{pmatrix}
-y^{(6)}+\sum\limits_{r=0}^2\left(g_ry^{(r)}\right)^{(r)}\\
U_{1}y\\
\end{pmatrix}
\text{ for } \widetilde{y} \in\mathscr{D}(A(U)), \\
K=\begin{pmatrix}
	0 & 0\\
	0 & I
\end{pmatrix}
\text{ and }
M=\begin{pmatrix}
I & 0\\
0 & 0
\end{pmatrix}.  
\end{gather*}
We can observe that $K\geq 0, M\geq 0$, $M+K=I$ and  $M|_{\mathscr{D}(A(U)}>0$.
 We associate a quadratic operator pencil
\begin{equation}\label{3.6}
L(\lambda) = \lambda^{2}M -i\lambda K-A(U), \quad \lambda\in\mathbb{C},
\end{equation}\\
\noindent in the space $L_{2}(a,b)\oplus\mathbb{C}^{l}$  with the problem \eqref{3.1}, \eqref{3.2}. The  operator pencil \eqref{3.6} is an operator realization of the eigenvalue problems \eqref{3.1}, \eqref{3.2} such that a function $y$ satisfies \eqref{3.1}, \eqref{3.2} if and only if $\widetilde y$ satisfies $L(\lambda)\widetilde{y}=0$.

We write the problem \eqref{3.1}, \eqref{3.2} in  the form
\begin{gather}\label{3.1A}
-y^{(6)}+\sum\limits_{r=0}^2\left(g_ry^{(r)}\right)^{(r)}=\lambda^{2}y,\\
\label{BC}
\qquad\quad  \left(\sum\limits_{k=1}^{6}\omega_{j,k}^{(0)}(\lambda)y^{(k-1)}(a)+\omega_{j,k}^{(1)}(\lambda)y^{(k-1)}(b)\right)_{j=1}^{6}=0.
\end{gather}
The representation of $\omega_{j,k}^{(0)}$ and $\omega_{j,k}^{(1)}$, $k=1,2,3,4,5,6,$, $j=1,2,3,4,5,6$ are
\begin{equation}\label{111}
\begin{cases}
\omega_{j,k}^{(0)}(\lambda)=\alpha_{j,k}+i\lambda\beta_{j,k}$  \text{ if }$j=1,2,3,
\\
\omega_{j,k}^{(0)}(\lambda)=0 \text{ if } j=4,5,6, 
\end{cases}
\end{equation}\\
\noindent   while

\begin{equation}\label{112}
\begin{cases}
\omega_{j,k}^{(1)}(\lambda)=0 \text{ if } j=1,2,3,\\
\omega_{j,k}^{(1)}(\lambda)=\alpha_{j,k}+i\lambda\beta_{j,k}$  \text{ if }$j=4,5,6,
\end{cases}
\end{equation}

\noindent where $k=1,\ldots,p_{j}$ for $\alpha_{j,k}$ and where $k=1,\ldots,q_{j}$ for $\beta_{j,k}$.

\section{Birkhoff regular problems}\label{Sec3}
 The characteristic function of the differential equation \eqref{3.1}, as defined in \cite[7.1.4]{mennicken2003non} is given by $\pi(\rho)=\rho^{6}+1$, and its roots are $\left(\sqrt{3}/{2}+i/{2}\right)^{(j-1)(2k+1)}, \text{ for } k,j=1,2,3,4,5,6.$

 We choose 
 \begin{equation}\label{C(x1)}
C(x,\mu)=\diag(1,\mu^{1},\mu^{2},\mu^{3},\mu^{4}, \mu^{5})\Big(\big(\sqrt{3}/{2}+{i}/{2})^{(j-1)(2k+1)}\Big),
\end{equation}
 see \cite[Theorem 7.2.4 A]{mennicken2003non}, where $\lambda =\mu^3$. The boundary matrices defined in \cite[(7.3.1)]{mennicken2003non} are

\begin{equation}
W^{(u)}=
(\omega_{j,k}^{(u)}(\mu^{3})C(a_u, \mu))_{j,k=1}^{6}, \quad u=0,1, 
\end{equation}
 where $a_u=0$ for $u=0$ and $a_u=a$ for $u=1$, are as defined in \eqref{111} and \eqref{112} for $\lambda=\mu^{3}$. Hence, 
\begin{equation}
\begin{aligned}
W^{(0)}(\mu) 
&=\begin{pmatrix}
\gamma_{1,k}\\
\gamma_{2,k}\\
\gamma_{3,k}\\
0\\
0\\
0\\
\end{pmatrix}_{k=1}^{6}, \quad W^{(1)}(\mu)=\begin{pmatrix}
0\\
0\\
0\\
\gamma_{4,k}\\
\gamma_{5,k}\\
\gamma_{6,k}\\
\end{pmatrix}_{k=1}^{6},
\end{aligned}
\end{equation}
 where
$\gamma_{j,k}=\sum\limits_{l=1}^{6}(\omega_{j,l}^{(u)})\xi^{(l-1)(2k+1)}\mu^{l-1},
$
\noindent with $\xi= \sqrt{3}/{2}+i/{2} $, $j=1,2,3$ for $u=0$ and $j=4,5,6$ for $u=1$.

Note that 
 $p_{j}$ and $q_{j}$ depends on the conditions stated in Assumption \ref{3.1.1}  and that $\beta_{j,q_{j}}\neq 0$ if $q_{j}\in \{0,1,2,3,4,5\}$, and $\alpha_{j,p_{j}}=1$ if $j\in\Lambda$, see (\ref{3.4}). 
 
Put 
 \begin{equation}\label{nu}
 \nu_{j}=\max\{p_{j},q_{j}+3\}.
 \end{equation}
 \noindent Choosing $C_{2}(\mu)=\text{diag}(\mu^{\nu_{1}},\mu^{\nu_{2}},\mu^{\nu_{3}}, \mu^{\nu_{4}}, \mu^{\nu_{5}}, \mu^{\nu_{6}})$ as defined in \cite{mennicken2003non}, we obtain  
  \begin{equation}
 C_{2}(\mu)^{-1}\cdot W^{(u)}(\mu) =W_{0}^{(u)}+O(\mu^{-1}), \quad u=0,1,
 \end{equation} 
 where 
 \begin{equation}\label{omega}
 W^{(0)}_{0}=\begin{pmatrix}
 \omega_{1,k}^{(0)}\\
 \omega_{2,k}^{(0)}\\
 \omega_{3,k}^{(0)}\\
  0\\
 0\\
 0\\
 \end{pmatrix}_{k=1}^{6}
 \text{ and } \ \ 
 W_{0}^{(1)}=\begin{pmatrix}
  0\\
 0\\
 0\\
 \omega_{4,k}^{(1)}\\
 \omega_{5,k}^{(1)}\\
 \omega_{6,k}^{(1)}\\
 \end{pmatrix}_{k=1}^{6}.
\end{equation}
 	The entries $\omega_{j,k}^{(u)}, u=0,1,$ represent the  coefficients of the terms with the highest degrees of polynomials $\gamma_{j ,k}$ in $\mu$, for $k=1,2,3,4,5,6$ and $j=1,2,3,4,5,6$.

 \noindent For simplicity, we set
 
 \begin{equation}\label{omega1}
 \omega_{j,k}^{(u)}=\omega_{j+3u,k+3u},
 \end{equation}\\
 for $k=1,2,3,4,5,6$ and $j=1,2,3,4,5,6$, with $u=0,1$. The Birkhoff matrices are given by
 
  \begin{equation}\label{birk}
W_{0}^{(0)}\Delta_{k}+W_{0}^{(1)}(I_{6}-\Delta_{k}), 
 \end{equation}
where $\Delta_k$, $k=1, \dots, 6$  are $6\times 6$ matrices with 3 consecutive ones and three consecutive zeros in the diagonal in a cyclic arrangement. After a permutation of columns, the matrices \eqref{birk} are blocks of diagonal matrices taken from three consecutive columns of the first three rows of $W_0^{(0)}$ and the last three rows of $W_0^{(1)}$ respectively. 
 
We define these block diagonal matrices  $\Gamma_{3u,k+3u}$\, 
 by
\begin{equation}\label{Gamma1}
\Gamma_{3u,k+3u}=
\begin{pmatrix}
\omega_{1+3u,k+3u} &\omega_{1+3u,k+1+3u} &\omega_{1+3u,k+2+3u}\\
\omega_{2+3u,k+3u} &\omega_{2+3u,k+1+3u} &\omega_{2+3u,k+2+3u}\\
\omega_{3+3u,k+3u} &\omega_{3+3u,k+1+3u} &\omega_{3+3u,k+2+3u}\\
\end{pmatrix},
\end{equation}
for $k=1,2,3,4,5,6$ and $u=0,1$. The indices of the entries of the matrix above can be simplified to $k+3u\equiv k+3u-6 \text{ mod } 6$ and $k+1+3u\equiv k+1+3u-6 \text{ mod } 6$, where $u=0,1$. 
The determinants of the Birkhoff matrices in \eqref{birk} can be written  in the form:

\begin{equation}
\det \left[W_{0}^{(0)}\Delta_{k}+W_{0}^{(1)}(I_{6}-\Delta_{k})\right]=\pm\det \Gamma_{0,k}\times\Gamma_{3,k+3}.\\
\end{equation}

 We can classify the eigenvalue problems \eqref{3.1}, \eqref{3.2} by the powers $p_{j}$ and $q_{j}$ of the derivatives in the boundary conditions \eqref{BC}. These classifications are given by:

\begin{equation}\label{classification1}
\left\{
\begin{array}{ll}     
\begin{aligned}
p_{j+3u}&<q_{j+3u}+3,\\
p_{j+3u}&>q_{j+3u}+3,\\
p_{j+3u}&=q_{j+3u}+3,\\
\end{aligned}
\end{array}\right.
\end{equation}
We have three different cases for each boundary conditions. Since each endpoint has three boundary conditions at each endpoint, then we have  27 distinct cases for each endpoint. From the above 27 distinct cases, three of the cases are not redundant, while  there are seven groups of redundant cases. The three non redundant  cases are  are given by:

\begin{itemize}
	\item[1.] $p_{1+3u}<q_{1+3u}+3,
	\qquad p_{2+3u}<q_{2+3u}+3,
	\qquad 	p_{3+3u}<q_{3+3u}+3$;
	\item[2.] $p_{1+3u}>q_{1+3u}+3,
	\qquad p_{2+3u}>q_{2+3u}+3,
	\qquad p_{3+3u}>q_{3+3u}+3$;
	\item[3.] $p_{1+3u}=q_{1+3u}+3,
	\qquad p_{2+3u}=q_{2+3u}+3,
	\qquad p_{3+3u}=q_{3+3u}+3$.
\end{itemize}

 Applying Assumption \ref{3.1.1} to the remaining 24 boundary conditions, we notice that 17 of them cannot hold. Hence, we are left with 7 additional boundary conditions.  These boundary conditions are given by

\begin{itemize}

	\item[1.] $p_{1+3u}>q_{1+3u}+3,
	\quad p_{2+3u}<q_{2+3u}+3,
	\quad p_{3+3u}<q_{3+3u}+3$;
	
		\item [2.] $p_{1+3u}>q_{1+3u}+3,
	\quad p_{2+3u}>q_{2+3u}+3,
	\quad p_{3+3u}<q_{3+3u}+3$;
	
\item[3.] $p_{1+3u}<q_{1+3u}+3,
	\quad p_{2+3u}<q_{2+3u}+3,
	\quad p_{3+3u}=q_{3+3u}+3$;

\item[4.] $p_{1+3u}>q_{1+3u}+3,
	\quad p_{2+3u}<q_{2+3u}+3,
	\quad p_{3+3u}=q_{3+3u}+3$;

	\item[5.]  $p_{1+3u}>q_{1+3u}+3,
	\quad p_{2+3u}>q_{2+3u}+3,
	\quad p_{3+3u}=q_{3+3u}+3$;

\item[6.] $p_{1+3u}>q_{1+3u}+3,
	\quad p_{2+3u}=q_{2+3u}+3,
	\quad p_{3+3u}=q_{3+3u}+3$;

	\item[7.] $p_{1+3u}<q_{1+3u}+3,	\quad 
	p_{2+3u}=q_{2+3u}+3,
	\quad p_{3+3u}=q_{3+3u}+3$.
\end{itemize}
 For practicality, we put
\begin{equation}\label{4.33}    
\begin{aligned}
&\theta_{j,u} =p_{j+3u}\quad \text{ and } \quad
 &\varphi_{j,u}=q_{j+3u}.\\
\end{aligned}
\end{equation}
Thus, the  conditions in \eqref{classification1} will be 
 
 \begin{equation}\label{classification2}
 \left\{
 \begin{array}{ll}     
 \begin{aligned}
 \theta_{j,u}&<\varphi_{j,u}+3,\\
  \theta_{j,u}&>\varphi_{j,u}+3,\\
  \theta_{j,u}&=\varphi_{j,u}+3.
 \end{aligned}
 \end{array}\right.
 \end{equation}
 We will denote the 10 different cases for each endpoint by Case$^{(u)}r$, where $u=0,1$ and $r=1,2, \dots,10$. Then we have
\begin{enumerate}
	\item 
	Case$^{(u)}$ 1: $\theta_{1,u}>\varphi_{1,u}+3,
	\quad \theta_{2,u}>\varphi_{2,u}+3,
	\quad \theta_{3,u}>\varphi_{3,u}+3$,
	
	\item Case$^{(u)}$ 2: $	\theta_{1,u}<\varphi_{1,u}+3,
	\quad 	\theta_{2,u}<\varphi_{2,u}+3,
	\quad \theta_{3,u}<\varphi_{3,u}+3$,
	
	\item Case$^{(u)}$ 3: $\theta_{1,u}>\varphi_{1,u}+3,
	\quad \theta_{2,u}>\varphi_{2,u}+3,
	\quad \theta_{3,u}<\varphi_{3,u}+3$,

	\item Case$^{(u)}$ 4: $\theta_{1,u}>\varphi_{1,u}+3,
	\quad \theta_{2,u}<\varphi_{2,u}+3,
	\quad \theta_{3,u}<\varphi_{3,u}+3$,

	\item Case$^{(u)}$ 5: $\theta_{1,u}>\varphi_{1,u}+3
	\quad \theta_{2,u}>\varphi_{2,u}+3,
	\quad \theta_{3,u}=\varphi_{3,u}+3$,

	\item Case$^{(u)}$ 6: $\theta_{1,u}>\varphi_{1,u}+3,
	\quad \theta_{2,u}<\varphi_{2,u}+3,
	\quad \theta_{3,u}=\varphi_{3,u}+3$,

	\item Case$^{(u)}$ 7: $\theta_{1,u}<\varphi_{1,u}+3,
	\quad \theta_{2,u}<\varphi_{2,u}+3,
	\quad \theta_{3,u}=\varphi_{3,u}+3$, 
	
	\item Case$^{(u)}$ 8: $\theta_{1,u}>\varphi_{1,u}+3,
	\quad \theta_{2,u}=\varphi_{2,u}+3,
	\quad \theta_{3,u}=\varphi_{3,u}+3$,
	
		\item Case$^{(u)}$ 9: $\theta_{1,u}<\varphi_{1,u}+3,
	\quad \theta_{2,u}=\varphi_{2,u}+3,
	\quad \theta_{3,u}=\varphi_{3,u}+3$,
	
	\item Case$^{(u)}$ 10: $\theta_{1,u}=\varphi_{1,u}+3,
	\quad \theta_{2,u}=\varphi_{2,u}+3,
	\quad \theta_{3,u}=\varphi_{3,u}+3$.
	
%
\end{enumerate}
\section{Classification of sixth order Birkhoff regular problems}\label{Sec4}
 We will now investigate $\det \Gamma_{3u,k+3u}$ as given in \eqref{Gamma1} for the 10 different cases. 
We use the following  conditions  
 
\noindent 1)  If $\theta_{j,u}>\varphi_{j,u}+3$ then we can deduce from
	$\gamma_{j ,k}$,	that the coefficients of the terms with the highest degrees of polynomials $\gamma_{j ,k}$ in $\mu$, for $k=1,2,3,4,5,6$, $j=1,2,3,4,5,6$,  and  $u=0,1$, is given by
\begin{equation}\label{4.34}
\begin{aligned}
\omega_{j+3u,k+3u}=\xi^{{\theta_{j,u}}(2(k+3u)+1)},
\end{aligned}
	\end{equation}
	 where $\theta_{j,u}\in\{0,1,2,3,4,5\}$ and $\alpha_{j,\theta_{j,u}}=1$.\\
\noindent 2) If $\theta_{j,u}<\varphi_{j,u}+3$, then we can deduce from
	$\gamma_{j ,k}$ that the coefficients of the terms with the highest degrees of polynomials $\gamma_{j ,k}$ in $\mu$, for $k=1,2,3,4,5,6$ and $j=1,2,3,4,5,6$, is given by
	\begin{equation}
\begin{aligned}\label{4.36}
	\omega_{j+3u,k+3u}&=i\xi^{{\varphi_{j,u}}(2(k+3u)+1)}\beta_{j+3u,\varphi_{j,u}},\\
\end{aligned}
	\end{equation}
 where $\varphi_{j,u}\in\{0,1,2,3,4,5\}$ and $\beta_{j,\varphi_{j,u}}\neq0$.\\
3)  If $	\theta_{j,u}=\varphi_{j,u}+3$, we can deduce from $\gamma_{j ,k}$ that the coefficients of the terms with the highest degrees of polynomials $\gamma_{j ,k}$ in $\mu$, for $k=1,2,3,4,5,6$, and $j=1,2,3,4,5,6$, is given by
	\begin{equation}\label{4.38}
	\omega_{j+3u,k+3u}=\xi^{{\varphi_{j,u}}(2(k+3u)+1)}({i^{2(k+3u)+1}+i\beta_{j+3u,\varphi_{j,u}}}),
	\end{equation}
		where $\theta_{j,u}, \varphi_{j,u}\in\{0,1,2,3,4,5\}$ with $\alpha_{j,\theta_{j,u}}=1$ and $\beta_{j,\varphi_{j,u}}\neq0$.
Putting
\begin{equation}\label{tau}
\tau_{k,u}=2(k+3u)+1,
\end{equation}                           
\noindent then

 \begin{equation}
\tau_{k+p,u}=\tau_{k,u}+2p, 
\end{equation} 
 \noindent  where $p$ is a natural number. Hence, 
\begin{lemma}\label{lem 4.4.28}\end{lemma}
\noindent 1) $i^{\tau_{k,u}}=\pm i$,\\
2)  $i^{\tau_{k,u}+1}=\pm 1$,\\
3) $i^{\tau_{k,u}+2}=\mp i$,\\
4)  $i^{2\tau_{k,u}}=-1$.

We will observe that
\begin{gather}\label{mzo1}\begin{cases}
\omega_{j+3u,k+3u}=\xi^{(\tau_{k,u})\theta_{j,u}},  \textrm{ for } \, \theta_{j,u}>\varphi_{j,u}+3,\\
\omega_{j+3u,k+3u}=i\xi^{\tau_{k,u}\varphi_{j,u}}\beta_{j+3u,\varphi_{j,u}}, \textrm{ for } \, \theta_{j,u}<\varphi_{j,u}+3,\\
\omega_{j+3u,k+3u}=\xi^{\tau_{k,u}\varphi_{j,u}}\phi_{j,\tau_{k,u}}, \textrm{ for }\, \theta_{j,u}=\varphi_{j,u}+3,
\end{cases}
\end{gather}
where
\begin{equation}\label{phi}
\phi_{j,\tau_{k,u}}=i^{\tau_{k,u}}+i\beta_{j+3u,\varphi_{j,u}}. 
\end{equation}
 
Using discrete Fourier transform, we write the matrices $\Gamma_{3u,k+3u}$ as
\begin{equation*}
\Gamma_{3u,k+3}=A\cdot \gamma_{3u,k+3u},
\end{equation*} 
where $A=\diag(\xi^mi\beta_{1+3u,\varphi_{1,u}}, \xi^ni\beta_{2+3u,\varphi_{2,u}}, \xi^qi\beta_{3+3u,\varphi_{3,u}})$ if $\theta_{j,u}<\varphi_{j,u}+3$. Observe that $i\beta_{j+3u,\varphi_{j,u}}=1$  if  $\theta_{j,u}\ge\varphi_{j,u}+3$, $j=1,2,3$ and $u=0,1$.
	It is easy to observe that $\xi^{\theta_{j,u}}, \, \xi^{\varphi_{j,u}}\ne 0$ for  $\theta_{j,u}, \varphi_{j,u},\in\{0,1,2,3,4,5\}$, $j=1,2,3,4,5,6$, $u=0,1$,   and $m, n, q\in\mathbb N$. Hence, it follows from \eqref{alphabeta} that $\det A\ne0$.

Since  $\det A\ne0$, comparing $\det \Gamma_{3u,k+3}$ to zero will be reduced to comparing $\det  \gamma_{3u,k+3}$ to zero. We will now evaluate $\det  \gamma_{3u,k+3}$, $u=0,1$, $k=1,2,3,4,5,6$.  For the above 10 cases, we obtain:

\noindent Case$^{(u)}$ 1: $\theta_{1,u}>\varphi_{1,u}+3,
	\quad \theta_{2,u}>\varphi_{2,u}+3,
	\quad \theta_{3,u}>\varphi_{3,u}+3$,
\begin{equation}\label{4.54}
\begin{aligned}
\gamma_{3u,k+3u}=&\begin{pmatrix}
1 &\xi^{2\theta_{1,u}} &\xi^{4\theta_{1,u}}\\
1 &\xi^{2\theta_{2,u}} &\xi^{4\theta_{2,u}}\\
1 &\xi^{2\theta_{3,u}} &\xi^{4\theta_{3,u}}
\end{pmatrix}
\end{aligned}
\end{equation}
and 
\begin{equation}\label{4.55}
\begin{aligned}
\det\gamma_{3u,k+3u}
=&(\xi^{2\theta_{2,u}}-\xi^{2\theta_{1,u}})(\xi^{2\theta_{3,u}}-\xi^{2\theta_{1,u}})(\xi^{2\theta_{3,u}}-\xi^{2\theta_{2,u}})\ne 0,
\end{aligned}
\end{equation}
 as $\theta_{1,u}, \theta_{2,u}$ and $\theta_{3,u}$, for $u$ $=0, 1$, are mutually disctinct, see Assumption \ref{3.1.1}. \\
Case$^{(u)}$ 2: $	\theta_{1,u}<\varphi_{1,u}+3,
	\quad 	\theta_{2,u}<\varphi_{2,u}+3,
	\quad \theta_{3,u}<\varphi_{3,u}+3$,
\begin{equation}
\gamma_{3u,k+3u}=\begin{pmatrix}
	1 &\xi^{2\varphi_{1,u}} &\xi^{4\varphi_{1,u}}\\
	1 &\xi^{2\varphi_{2,u}} &\xi^{4\varphi_{2,u}}\\
	1 &\xi^{2\varphi_{3,u}} &\xi^{4\varphi_{3,u}}
\end{pmatrix}
\end{equation}
and it follows from  \eqref{alphabeta}
\begin{equation}\label{4.63}
\begin{aligned}
\det\gamma_{3u,k+3u}
=&(\xi^{2\varphi_{2,u}}-\xi^{2\varphi_{1,u}})(\xi^{2\varphi_{3,u}}-\xi^{2\varphi_{1,u}})(\xi^{2\varphi_{3,u}}-\xi^{2\varphi_{2,u}})\ne0,
\end{aligned}
\end{equation}
since $\varphi_{j,u}$, $j=1,2,3$ nd $u=0,1$, are mutually distinct.\\
Case$^{(u)}$ 3: $\theta_{1,u}>\varphi_{1,u}+3,
	\quad \theta_{2,u}>\varphi_{2,u}+3,
	\quad \theta_{3,u}<\varphi_{3,u}+3$,
\begin{equation}
\begin{aligned}
\gamma_{3u,k+3u}&=
\begin{pmatrix}
1 &\xi^{2\theta_{1,u}} &\xi^{4\theta_{1,u}}\\
1 &\xi^{2\theta_{2,u}} &\xi^{4\theta_{2,u}}\\
1 &\xi^{2\varphi_{3,u}} &\xi^{4\varphi_{3,u}}
\end{pmatrix}.
\end{aligned}
\end{equation}
Thus,
\begin{equation}\label{4.691}
\begin{aligned}
\det\gamma_{3u,k+3u}=
&(\xi^{2\varphi_{2,u}}-\xi^{2\theta_{1,u}})(\xi^{2\varphi_{3,u}}-\xi^{2\theta_{1,u}})(\xi^{2\theta_{3,u}}-\xi^{2\varphi_{2,u}})\ne0.
\end{aligned}
\end{equation}
Case$^{(u)}$ 4: $\theta_{1,u}>\varphi_{1,u}+3,
	\quad \theta_{2,u}<\varphi_{2,u}+3,
	\quad \theta_{3,u}<\varphi_{3,u}+3$,

\begin{equation}
\gamma_{3u,k+3u}=\begin{pmatrix}
1 &\xi^{2\theta_{u}} &\xi^{4\theta_{u}}\\
1 &\xi^{2\varphi_{2u}} &\xi^{4\varphi_{2 u}}\\
1 &\xi^{2\varphi_{3u}} &\xi^{4\varphi_{3u}}
\end{pmatrix}.
\end{equation}
Since $\theta_{j,u}$, $\varphi_{k,u}$, $j,k=1,2,3$ and $u=0,1$, are mutually disctinct, 
\begin{equation}\label{4.74}
\begin{aligned}
\det\gamma_{3u,k+3u}
&=(\xi^{2\varphi_{2,u}}-\xi^{2\theta_{1,u}})(\xi^{2\varphi_{3,u}}-\xi^{2\theta_{1,u}})(\xi^{2\varphi_{2,u}}-\xi^{2\varphi_{3,u}})\ne0.
\end{aligned}
\end{equation}
From the above analyses,  the following result follows:
\begin{proposition}\label{4.5.1}
	Let $\theta_{j,u}\in\{-\infty,0,1,2,3,4,5\}$ and $\varphi_{j,u}\in \{-\infty, 0,1,2\}$, where $\theta_{j,u}=p_{j+3u}$ and $\varphi_{j,u}=q_{j+3u}$ are mutually exclusive while  $j=1,2,3$ and  $u=0,1$.  Then
	\begin{equation*}
	\det\Gamma_{3u,k+3u}\neq 0
	\end{equation*}
	\noindent \textit{for the following conditions:}
	\begin{enumerate}
		\item $\theta_{1,u}>\varphi_{1,u}+3,
		\quad \theta_{2,u}>\varphi_{2,u}+3,
		\quad \theta_{3,u}>\varphi_{3,u}+3$,
		\item $\theta_{1,u}<\varphi_{1,u}+3,
		\quad \theta_{2,u}<\varphi_{2,u}+3,
		\quad \theta_{3,u}<\varphi_{3,u}+3$,
		\item $\theta_{1,u}>\varphi_{1,u}+3,
		\quad \theta_{2,u}>\varphi_{2,u}+3,
		\quad \theta_{3,u}<\varphi_{3,u}+3$,
		\item 	$\theta_{1,u}>\varphi_{1,u}+3,
		\quad \theta_{2,u}<\varphi_{2,u}+3,
		\quad \theta_{3,u}<\varphi_{3,u}+3$.
	\end{enumerate}
\end{proposition}
Next, we investigate the remaining  cases where conditions are needed for \\$\det\Gamma_{2u,k+2u}\neq 0$.

\noindent Case$^{(u)}$ 5: $\theta_{1,u}>\varphi_{1,u}+3
	\quad \theta_{2,u}>\varphi_{2,u}+3,
	\quad \theta_{3,u}=\varphi_{3,u}+3$,
\begin{equation}\label{gamma5}
\gamma_{3u,k+3u}=\begin{pmatrix}
	1 &\xi^{2\theta_{1,u}} &\xi^{4\theta_{1,u}}\\
	1 &\xi^{2\theta_{2,u}} &\xi^{4\theta_{2,u}}\\
	\phi_{3,\tau_{k,u}} &\xi^{2\varphi_{3,u}}\phi_{3,\tau_{k,u}+2} &\xi^{4\varphi_{3,u}}\phi_{3,\tau_{k,u}}
\end{pmatrix}.
\end{equation}
Observe from Lemma \ref{lem 4.4.28} that 
	\begin{equation}\label{4.69}
	\phi_{j,\tau_{k,u}}=0 \Longleftrightarrow  \beta_{3+3u,\varphi_{j,u}}=\pm 1, ~u=0,1.
	\end{equation}
	\noindent  Hence,
\begin{equation}\label{5.1}
\begin{aligned}
\det\gamma_{3u,k+3u}&=0.
\end{aligned}
\end{equation}
For  $\beta_{3+3u, \varphi_{3,u}}\ne\pm 1$, 
\begin{align}\label{4.67}
\det\gamma_{3u,k+3u}=&(\xi^{2\theta_{2,u}}-\xi^{2\theta_{1,u}})(\xi^{4\varphi_3,u}\phi_{3,\tau_{k,u}}-\xi^{4\theta_{1,u}}\phi_{3,\tau_{k,u}})\\\nonumber
&-(\xi^{4\theta_{2,u}}-\xi^{4\theta_{1,u}})(\xi^{2\varphi_3,u}\phi_{3,\tau_{k,u}+2}-\xi^{2\theta_{1,u}}\phi_{3,\tau_{k,u}})\nonumber\\
&=(\xi^{2\theta_{2,u}}-\xi^{2\theta_{1,u}})(\xi^{2\varphi_{3,u}}+\xi^{2\theta_{2,u}})\nonumber\\&\qquad \times \Big[i^{\tau_{3,u}}(\xi^{2\varphi_{3,u}}+\xi^{2\theta_{1,u}})+i\beta_{3+3u, \varphi_{3,u}}(\xi^{2\varphi_{3,u}}-\xi^{2\theta_{2,u}})\Big].
\end{align}
Since $\theta_1$, $\theta_2$ and $\varphi_3$ are mutually distinct and $\varphi_3\in\{0,1,2\}$, then 
\begin{equation}\label{moz51}\det\gamma_{3u,k+3u}=0 \iff \beta_{3+3u, \varphi_{3,u}}=\mp\dfrac{\xi^{2\varphi_3}+\xi^{2\theta_1} }{\xi^{2\varphi_3}-\xi^{2\theta_2}}.\end{equation} Hence,
\begin{proposition}\label{prop5}
	Let $\xi= \sqrt{3}/{2}+i/{2} $, $\theta_{j,u}\in\{-\infty,0,1,2,3,4,5\}$ and $\varphi_{j,u}\in \{-\infty, 0,1,2\}$,  where $\theta_{j,u}=p_{j+3u}$ and $\varphi_{j,u}=q_{j+3u}$ are mutually exclusive  while  $j=1,2,3$, $u=0,1$.   \\
 For $\theta_{1,u}>\varphi_{1,u}+3, \, 
	 \theta_{2,u}>\varphi_{2,u}+3, \,
	 \theta_{3,u}=\varphi_{3,u}+3$,
\begin{gather*}
\det\Gamma_{3u,k+3u}\neq 0 \Longleftrightarrow 
\beta_{3+3u,\varphi_{j,u}}\ne
\begin{cases}
\pm1,\\
\mp\dfrac{\xi^{2\varphi_3}+\xi^{2\theta_1} }{\xi^{2\varphi_3}-\xi^{2\theta_2}}.
\end{cases}
\end{gather*}
\end{proposition}
\noindent Case$^{(u)}$ 6: $\theta_{1,u}>\varphi_{1,u}+3,
	\quad \theta_{2,u}<\varphi_{2,u}+3,
	\quad \theta_{3,u}=\varphi_{3,u}+3$,
\begin{equation}\label{gamma6}
\gamma_{3u,k+3u}=\begin{pmatrix}
1 &\xi^{2\theta_{1,u}} &\xi^{4\theta_{1,u}}\\
1 &\xi^{2\varphi_{2,u}} &\xi^{4\varphi_{2,u}}\\
\phi_{3,\tau_{k,u}} &\xi^{2\varphi_{3,u}}\phi_{3,\tau_{k,u}+2} &\xi^{4\varphi_{3,u}}\phi_{3,\tau_{k,u}}\\
\end{pmatrix}.
\end{equation}	
It follows from \eqref{gamma5}, \eqref{gamma6} and Proposition \ref{prop5} the following result
\begin{proposition}\label{prop6}
	Let $\xi= \sqrt{3}/{2}+i/{2} $, $\theta_{j,u}\in\{-\infty,0,1,2,3,4,5\}$ and $\varphi_{j,u}\in \{-\infty, 0,1,2\}$,  where $\theta_{j,u}=p_{j+3u}$ and $\varphi_{j,u}=q_{j+3u}$ are mutually exclusive  while  $j=1,2,3$,  $u=0,1$.  \\For  Case$^{(u)}$ 6: $\theta_{1,u}>\varphi_{1,u}+3,
	\quad \theta_{2,u}<\varphi_{2,u}+3,
	\quad \theta_{3,u}=\varphi_{3,u}+3$,
	\begin{gather*}
\det\Gamma_{3u,k+3u}\neq 0 \Longleftrightarrow 
\beta_{3+3u,\varphi_{j,u}}\ne
\begin{cases}
\pm1,\\
\mp\dfrac{\xi^{2\varphi_3}+\xi^{2\theta_1} }{\xi^{2\varphi_3}-\xi^{2\varphi_2}}.
\end{cases}
\end{gather*}
\end{proposition}

\noindent Case$^{(u)}$ 7: $\theta_{1,u}<\varphi_{1,u}+3,
	\quad \theta_{2,u}<\varphi_{2,u}+3,
	\quad \theta_{3,u}=\varphi_{3,u}+3$, 
\begin{equation}\label{gamma7}
\gamma_{3u,k+3u}=\begin{pmatrix}
1 &\xi^{2\varphi_{1,u}} &\xi^{4\varphi_{1,u}}\\
1 &\xi^{2\varphi_{2,u}} &\xi^{4\varphi_{2,u}}\\
\phi_{3,\tau_{k,u}} &\xi^{2\varphi_{3,u}}\phi_{3,\tau_{k,u}+2} &\xi^{4\varphi_{3,u}}\phi_{3,\tau_{k,u}}
\end{pmatrix}.
\end{equation}	
Hence, from \eqref{gamma5} and Proposition \ref{prop5}, we have the following result,
\begin{proposition}\label{prop7}
	Let $\xi= \sqrt{3}/{2}+i/{2} $, $\theta_{j,u}\in\{-\infty,0,1,2,3,4,5\}$ and $\varphi_{j,u}\in \{-\infty, 0,1,2\}$,  where $\theta_{j,u}=p_{j+3u}$ and $\varphi_{j,u}=q_{j+3u}$ are mutually exclusive while  $j=1,2,3$,  $u=0,1$.   \\For  Case$^{(u)}$ 7: $\theta_{1,u}<\varphi_{1,u}+3,
	\quad \theta_{2,u}<\varphi_{2,u}+3,
	\quad \theta_{3,u}=\varphi_{3,u}+3$, 
\begin{gather*}
\det\Gamma_{3u,k+3u}\neq 0 \Longleftrightarrow 
\beta_{3+3u,\varphi_{j,u}}\ne
\begin{cases}
\pm1,\\
\mp\dfrac{\xi^{2\varphi_3}+\xi^{2\varphi_1} }{\xi^{2\varphi_3}-\xi^{2\varphi_2}}.
\end{cases}
\end{gather*}
	\end{proposition}
\noindent Case$^{(u)}$ 8: $\theta_{1,u}>\varphi_{1,u}+3,
	\quad \theta_{2,u}=\varphi_{2,u}+3,
	\quad \theta_{3,u}=\varphi_{3,u}+3$,
	\begin{equation}\label{det9}
\begin{aligned}
\gamma_{3u,k+3u}=
\begin{pmatrix}
1 &\xi^{2\theta_{1,u}} &\xi^{4\theta_{1,u}}\\
\phi_{2,\tau_{k,u}} &\xi^{2\varphi_{2,u}}\phi_{2,\tau_{k,u}+2} &\xi^{4\varphi_{2,u}}\phi_{2,\tau_{k,u}}\\
\phi_{3,\tau_{k,u}} &\xi^{2\varphi_{3,u}}\phi_{3,\tau_{k,u}+2} &\xi^{4\varphi_{3,u}}\phi_{3,\tau_{k,u}}
\end{pmatrix}.
\end{aligned}  
\end{equation} 
Hence,
\begin{align}\label{Case8}
\det\gamma_{3u,k+3u}=&(\xi^{2\varphi_{2,u}}\phi_{2,\tau_{k,u}+2}-\xi^{2\theta_{1,u}}\phi_{2,\tau_{k,u}})(\xi^{4\varphi_3,u}\phi_{3,\tau_{k,u}}-\xi^{4\theta_{1,u}}\phi_{3,\tau_{k,u}})\\\nonumber
&-(\xi^{4\varphi_{2,u}}\phi_{2,\tau_{k,u}}-\xi^{4\theta_{1,u}}\phi_{2,\tau_{k,u}})(\xi^{2\varphi_3,u}\phi_{3,\tau_{k,u}+2}-\xi^{2\theta_{1,u}}\phi_{3,\tau_{k,u}}).
\end{align}
From  $\beta_{j+3u, \varphi_{3,u}}\ne\pm 1$, $j=2,3$, and  Assumption \ref{3.1.1}, we observe the following subcases

\noindent 1) Subcase$^{(u)}$  8.1: $\theta_1=2$, $\varphi_2=0$, $\varphi_3=1$,\\
2) Subcase$^{(u)}$  8.2: $\theta_1=5$, $\varphi_2=0$, $\varphi_3=1$,\\
3) Subcase$^{(u)}$  8.3: $\theta_1=1$, $\varphi_2=0$, $\varphi_3=2$,\\
4) Subcase$^{(u)}$  8.4: $\theta_1=4$, $\varphi_2=0$, $\varphi_3=2$,\\
5) Subcase$^{(u)}$  8.5: $\theta_1=0$, $\varphi_2=1$, $\varphi_3=2$,\\
6) Subcase$^{(u)}$  8.6: $\theta_1=3$, $\varphi_2=1$, $\varphi_3=2$.

We next investigate the different subcases from Case 8$^{(u)}$.\\
\noindent 1) {\bf Subcase$^{(u)}$ 8.1}:   $\theta_1=2$, $\varphi_2=0$, $\varphi_3=1$,
	\begin{align}\label{9mzeq1}\det\gamma_{3u,k+3u}&=-\sqrt3i\beta_{3+3u, \varphi_{3,u}}(i^{\tau_{k,u}+1}+\beta_{2+3u, \varphi_{2,u}})\nonumber\\&\quad +((-1)^{\tau_{k,u}}+3i^{\tau_{k,u}+1}\beta_2)\xi^4+((-1)^{\tau_{k,u}}+2i^{\tau_{k,u}+1})\xi^2.\end{align}
It follows from Lemma \ref{lem 4.4.28}, that
\begin{align}\label{subcase8.1}\det\gamma_{3u,k+3u}=0&\Longleftrightarrow\nonumber\\& \,~\beta_{3+3u, \varphi_{3,u}}=\dfrac{((-1)^{\tau_{k,u}}+3i^{\tau_{k,u}+1}\beta_2)\xi^4+((-1)^{\tau_{k,u}}+2i^{\tau_{k,u}+1})\xi^2}{\sqrt3i(i^{\tau_{k,u}+1}+\beta_{2+3u, \varphi_{2,u}})}.\end{align}
\begin{lemma}\label{lem8.1}
Let $\xi= \sqrt{3}/{2}+i/{2} $, $\theta_{j,u}\in\{-\infty,0,1,2,3,4,5\}$ and $\varphi_{j,u}\in \{-\infty, 0,1,2\}$,  where $\theta_{j,u}=p_{j+3u}$ and $\varphi_{j,u}=q_{j+3u}$ are mutually exclusive while  $j=1,2,3$, $u=0,1$.  Let $$ 
\beta_{j+3u,\varphi_{j,u}}\ne \pm1, \, j=2,3.$$
For  $\theta_1=2$, $\varphi_2=0$, $\varphi_3=1$,
\begin{align}\label{lemSubcase8.1}\det\Gamma_{3u,k+3u}\ne0&\Longleftrightarrow\nonumber\\&\quad  \beta_{3+3u, \varphi_{3,u}}\ne\dfrac{(-1)^{\tau_{k,u}}+3i^{\tau_{k,u}+1}\beta_2)\xi^4+((-1)^{\tau_{k,u}}+2i^{\tau_{k,u}+1})\xi^2}{\sqrt3i(i^{\tau_{k,u}+1}+\beta_{2+3u, \varphi_{2,u}})}.\end{align}
\end{lemma}
\noindent 2) {\bf Subcase$^{(u)}$ 8.2}:  $\theta_1=5$, $\varphi_2=0$, $\varphi_3=1$, 
\begin{align}\label{det8.2}
\det \gamma_{3u,k+3u}	=& \sqrt3i\Big(\beta_{3+3u, \varphi_{3,u}}(\beta_{2+3u, \varphi_{2,u}}-3i^{\tau_{k,u}+1})\nonumber\\&\qquad+((-1)^{\tau_{k,u}+1}+i^{\tau_{k,u}+1}\beta_{2+3u, \varphi_{2,u}})\Big).
\end{align}
If $\beta_{2+3u, \varphi_{2,u}}= \pm3$, then $\det \gamma_{3u,k+3u}\ne 0$. 
    However, if  $\beta_{2+3u, \varphi_{2,u}}\ne \pm3$, then 
\begin{align}\label{det8.2.1}\det \gamma_{3u,k+3u}\ne 0\Longleftrightarrow \beta_{3+3u, \varphi_{3,u}}\ne\dfrac{(-1)^{\tau_{k,u}+1}+i^{\tau_{k,u}+1}\beta_{2+3u, \varphi_{2,u}}}{3i^{\tau_{k,u}+1}-\beta_{2+3u, \varphi_{2,u}}}.\end{align}
\begin{lemma}\label{lem8.2}
Let $\xi= \sqrt{3}/{2}+i/{2} $, $\theta_{j,u}\in\{-\infty,0,1,2,3,4,5\}$ and $\varphi_{j,u}\in \{-\infty, 0,1,2\}$,  where $\theta_{j,u}=p_{j+3u}$ and $\varphi_{j,u}=q_{j+3u}$ are mutually exclusive while  $j=1,2,3$,  $u=0,1$.  Let  
$$ 
\beta_{j+3u,\varphi_{j,u}}\ne \pm1, \, j=2,3.$$
 For $\theta_1=5, \,\varphi_2=0, \, \varphi_3=1,$
 $$\det\Gamma_{3u,k+3u}\neq 0$$
 if and only one of the following conditions holds:\\
 1) $\beta_{2+3u,\varphi_{2,u}}=\pm 3$,\\
 2) $\beta_{2+3u,\varphi_{2,u}}\ne\pm 3$ and $\beta_{3+3u,\varphi_{3,u}}\ne\dfrac{(-1)^{\tau_{k,u}+1}+i^{\tau_{k,u}+1}\beta_{2+3u, \varphi_{2,u}}}{3i^{\tau_{k,u}+1}-\beta_{2+3u, \varphi_{2,u}}}$.
\end{lemma}
\noindent 3) {\bf Subcase$^{(u)}$  8.3}: $\theta_1=1$, $\varphi_2=0$, $\varphi_3=2$,
\begin{align}\label{det8.3}
\det \gamma_{3u,k+3u}	=& \Big(\beta_{3+3u, \varphi_{3,u}}(\beta_{2+3u, \varphi_{2,u}}+i^{\tau_{k,u}+1})\nonumber\\&\qquad+(3(-1)^{\tau_{k,u}+1}+i^{\tau_{k,u}+1}\beta_{2+3u, \varphi_{2,u}})\Big)(\xi^4+\xi^2).
\end{align}
Since $\beta_{2+3u, \varphi_{2,u}}\ne \pm1$, then 
\begin{align}\label{det8.3.1}\det \gamma_{3u,k+3u}\ne 0\Longleftrightarrow \beta_{3+3u, \varphi_{3,u}}\ne-\dfrac{3(-1)^{\tau_{k,u}+1}+i^{\tau_{k,u}+1}\beta_{2+3u, \varphi_{2,u}}}{i^{\tau_{k,u}+1}+\beta_{2+3u, \varphi_{2,u}}}.\end{align}
\begin{lemma}\label{lem8.3}
Let $\xi= \sqrt{3}/{2}+i/{2} $, $\theta_{j,u}\in\{-\infty,0,1,2,3,4,5\}$ and $\varphi_{j,u}\in \{-\infty, 0,1,2\}$,  where $\theta_{j,u}=p_{j+3u}$ and $\varphi_{j,u}=q_{j+3u}$ are mutually exclusive while  $j=1,2,3$,  $u=0,1$.  Let $$ 
\beta_{j+3u,\varphi_{j,u}}\ne \pm1, \, j=2,3.$$
For
$\theta_1=1, \, \varphi_2=0, \,\varphi_3=2.$ Then
\begin{align*}\det \Gamma_{3u,k+3u}\ne 0\Longleftrightarrow \beta_{3+3u, \varphi_{3,u}}\ne-\dfrac{3(-1)^{\tau_{k,u}+1}+i^{\tau_{k,u}+1}\beta_{2+3u, \varphi_{2,u}}}{i^{\tau_{k,u}+1}+\beta_{2+3u, \varphi_{2,u}}}.\end{align*}
\end{lemma}
\noindent 4) {\bf Subcase$^{(u)}$  8.4}: $\theta_1=4$, $\varphi_2=0$, $\varphi_3=2$,
\begin{align}\label{det8.4}
\det \gamma_{3u,k+3u}	=&  \Big(\beta_{3+3u, \varphi_{3,u}}(\beta_{2+3u, \varphi_{2,u}}-3i^{\tau_{k,u}+1})\xi^4+(i^{\tau_{k,u}+1}+\beta_{2+3u, \varphi_{2,u}})\xi^2\nonumber\\&\qquad+(i^{\tau_{k,u}+1}\beta_{2+3u, \varphi_{2,u}}+(-1)^{\tau_{k,u}+1})\xi^4\nonumber\\&\qquad\qquad+((-1)^{\tau_{k,u}+1}-3i^{\tau_{k,u}+1}\beta_{2+3u, \varphi_{2,u}})\xi^2.
\end{align}
If $\beta_{2+3u, \varphi_{2,u}}= \pm3$, then $\det \gamma_{3u,k+3u}\ne 0$. 
    However, if  $\beta_{2+3u, \varphi_{2,u}}\ne \pm3$, then 
    \begin{align}\label{det8.4.1}\det \gamma_{3u,k+3u}\ne 0\Longleftrightarrow \beta_{3+3u, \varphi_{3,u}}\ne&-\dfrac{(i^{\tau_{k,u}+1}+\beta_{2+3u, \varphi_{2,u}})(\xi^2+i^{\tau_{k,u}+1}\xi^4)}{3i^{\tau_{k,u}+1}-\beta_{2+3u, \varphi_{2,u}}}\nonumber\\&\qquad-\dfrac{((-1)^{\tau_{k,u}+1}-3i^{\tau_{k,u}+1}\beta_{2+3u, \varphi_{2,u}})\xi^2}{3i^{\tau_{k,u}+1}-\beta_{2+3u, \varphi_{2,u}}}.\end{align}
    \begin{lemma}\label{lem8.4}
Let $\xi= \sqrt{3}/{2}+i/{2} $, $\theta_{j,u}\in\{-\infty,0,1,2,3,4,5\}$ and $\varphi_{j,u}\in \{-\infty, 0,1,2\}$,  where $\theta_{j,u}=p_{j+3u}$ and $\varphi_{j,u}=q_{j+3u}$ are mutually exclusive while  $j=1,2,3$,  $u=0,1$.  Let $$ 
\beta_{j+3u,\varphi_{j,u}}\ne \pm1, \, j=2,3.$$
and 
    \begin{align*} \alpha_{3+3u, \varphi_{3,u}}=&-\dfrac{(i^{\tau_{k,u}+1}+\beta_{2+3u, \varphi_{2,u}})(\xi^2+i^{\tau_{k,u}+1}\xi^4)}{3i^{\tau_{k,u}+1}-\beta_{2+3u, \varphi_{2,u}}}\nonumber\\&\qquad-\dfrac{((-1)^{\tau_{k,u}+1}-3i^{\tau_{k,u}+1}\beta_{2+3u, \varphi_{2,u}})\xi^2}{3i^{\tau_{k,u}+1}-\beta_{2+3u, \varphi_{2,u}}}.\end{align*}
For $\theta_1=4$, $\varphi_2=0$, $\varphi_3=2$,
$$\det\Gamma_{3u,k+3u}\neq 0 $$
if and only one of the following conditions holds:\\
1) $\beta_{2+3u,\varphi_{2,u}}=\pm 3$,\\
2) $\beta_{2+3u,\varphi_{2,u}}\ne\pm 3$ and $\beta_{3+3u,\varphi_{3,u}}\ne\alpha_{3+3u, \varphi_{3,u}}$.
\end{lemma}
\noindent     5) {\bf Subcase$^{(u)}$  8.5}: $\theta_1=0$, $\varphi_2=1$, $\varphi_3=2$,
      \begin{align}\label{det8.5.1}	\det\gamma_{3u,k+3u}&=\beta_{3+3u,\varphi_{3,u}}(-\beta_{2+3u,\varphi_{2,u}}(\xi^4+\xi^2)+3i^{\tau_{k,u}+1}\xi^4-i^{\tau_{k,u}+1}\xi^2)\nonumber\\&\qquad+((-1)^{\tau_{k,u}+1}-i^{\tau_{k,u}+1}\beta_{2+3u,\varphi_{2,u}})\xi^4\nonumber\\&\qquad-3((-1)^{\tau_{k,u}+1}+\beta_{2+3u,\varphi_{2,u}})\xi^2
  \end{align}
  If $\beta_{2+3u, \varphi_{2,u}}=\dfrac{3i^{\tau_{k,u}+1}\xi^4-i^{\tau_{k,u}+1}\xi^2}{\sqrt3i}$, then $\det \gamma_{3u,\varphi_{k+3u}}\ne0$. However, if \\$\beta_{2+3u, \varphi_{2,u}}\ne\dfrac{3i^{\tau_{k,u}+1}\xi^4-i^{\tau_{k,u}+1}\xi^2}{\sqrt3i}$, then 
  \begin{align}\label{det8.5}&\det \gamma_{3u,\varphi_{k+3u}}\ne0\Longleftrightarrow \nonumber\\~&\beta_{3+3u, \varphi_{3,u}}\ne \dfrac{3((-1)^{\tau_{k,u}+1}+\beta_{2+3u,\varphi_{j,u}})\xi^2+(i^{\tau_{k,u}+1}\beta_{2+3u,\varphi_{j,u}}-(-1)^{\tau_{k,u}+1})\xi^4}{3i^{\tau_{k,u}+1}\xi^4-i^{\tau_{k,u}+1}\xi^2+\sqrt3i\beta_{2+3u, \varphi_{2,u}}}.\end{align}
   \begin{lemma}\label{lem8.5}
Let $\xi= \sqrt{3}/{2}+i/{2} $, $\theta_{j,u}\in\{-\infty,0,1,2,3,4,5\}$ and $\varphi_{j,u}\in \{-\infty, 0,1,2\}$,  where $\theta_{j,u}=p_{j+3u}$ and $\varphi_{j,u}=q_{j+3u}$ are mutually exclusive while  $j=1,2,3$,  $u=0,1$.  Let $$ 
\beta_{j+3u,\varphi_{j,u}}\ne \pm1, \, j=2,3.$$
and 
    \begin{align*} \alpha_{3+3u, \varphi_{3,u}}=&\dfrac{3((-1)^{\tau_{k,u}+1}+\beta_{2+3u,\varphi_{2,u}})\xi^2+(i^{\tau_{k,u}+1}\beta_{2+3u,\varphi_{j,u}}-(-1)^{\tau_{k,u}+1})\xi^4}{3i^{\tau_{k,u}+1}\xi^4-i^{\tau_{k,u}+1}\xi^2+\sqrt3i\beta_{2+3u, \varphi_{2,u}}}.\end{align*}
For  $\theta_1=0$, $\varphi_2=1$, $\varphi_3=2$,
$$\det\Gamma_{3u,k+3u}\neq 0 $$
if and only one of the following conditions holds:\\\\
1) $\beta_{2+3u,\varphi_{j,u}}=\dfrac{3i^{\tau_{k,u}+1}\xi^4-i^{\tau_{k,u}+1}\xi^2}{\sqrt3i}$,\\\\
2) $\beta_{2+3u,\varphi_{j,u}}\ne\dfrac{3i^{\tau_{k,u}+1}\xi^4-i^{\tau_{k,u}+1}\xi^2}{\sqrt3i}$ and $\beta_{3+3u,\varphi_{j,u}}\ne\alpha_{3+3u, \varphi_{3,u}}$.

\end{lemma}
 \noindent   6) {\bf Subcase$^{(u)}$  8.6}: $\theta_1=3$, $\varphi_2=1$, $\varphi_3=2$,
      \begin{align}\label{det8.6.1}	\det\gamma_{3u,k+3u}&=\beta_{3+3u,\varphi_{j,u}}(2i^{\tau_{k,u}+1}\xi^4+2\beta_{2+3u,\varphi_{j,u}}\xi^2)-2i^{\tau_{k,u}+1}\beta_{2+3u,\varphi_{j,u}}\xi^4\nonumber\\&\qquad -(4i^{\tau_{k,u}+1}\beta_{2+3u,\varphi_{j,u}}+(-1)^{\tau_{k,u}+1})\xi^2.
    \end{align}
    If $\beta_{2+3u, \varphi_{2,u}}=-i^{\tau_{k,u}+1}\xi^2$, then $\det \gamma_{3u,\varphi_{k+3u}}\ne0$.  On the other hand if $\beta_{2+3u, \varphi_{2,u}}\ne-i^{\tau_{k,u}+1}\xi^2$, then 
  \begin{align}\label{det8.3.1}&\det \gamma_{3u,k+3u}\ne 0\Longleftrightarrow\nonumber\\& \quad\beta_{3+3u, \varphi_{3,u}}\ne\dfrac{2i^{\tau_{k,u}+1}\beta_{2+3u,\varphi_{2,u}}\xi^4+(4i^{\tau_{k,u}+1}\beta_{2+3u,\varphi_{2,u}}+(-1)^{\tau_{k,u}+1})\xi^2}{2i^{\tau_{k,u}+1}\xi^4+2\beta_{2+3u,\varphi_{2,u}}\xi^2}.\end{align}
   \begin{lemma}\label{lem8.6}
Let $\xi= \sqrt{3}/{2}+i/{2} $, $\theta_{j,u}\in\{-\infty,0,1,2,3,4,5\}$ and $\varphi_{j,u}\in \{-\infty, 0,1,2\}$,  where $\theta_{j,u}=p_{j+3u}$ and $\varphi_{j,u}=q_{j+3u}$ are mutually exclusive while  $j=1,2,3$,   $u=0,1$.  Let $$ 
\beta_{j+3u,\varphi_{j,u}}\ne \pm1, \, j=2,3.$$
and 
    \begin{align*} \alpha_{3+3u, \varphi_{3,u}}=&\dfrac{2i^{\tau_{k,u}+1}\beta_{2+3u,\varphi_{2,u}}\xi^4+(4i^{\tau_{k,u}+1}\beta_{2+3u,\varphi_{2,u}}+(-1)^{\tau_{k,u}+1})\xi^2}{2i^{\tau_{k,u}+1}\xi^4+2\beta_{2+3u,\varphi_{2,u}}\xi^2}.\end{align*}
For $\theta_1=3$, $\varphi_2=1$, $\varphi_3=2$,
$$\det\Gamma_{3u,k+3u}\neq 0 $$
if and only one of the following conditions holds:\\
1) $\beta_{2+3u,\varphi_{2,u}}=-i^{\tau_{k,u}+1}\xi^2$,\\
2) $\beta_{2+3u,\varphi_{2,u}}\ne-i^{\tau_{k,u}+1}\xi^2$ and $\beta_{3+3u,\varphi_{3,u}}\ne\alpha_{3+3u, \varphi_{3,u}}$.

\end{lemma}
  \begin{proposition}\label{prop8}
  	Let $\xi= \sqrt{3}/{2}+i/{2} $, $\theta_{j,u}\in\{-\infty,0,1,2,3,4,5\}$ and $\varphi_{j,u}\in \{-\infty, 0,1,2\}$,  where $\theta_{j,u}=p_{j+3u}$ and $\varphi_{j,u}=q_{j+3u}$ are mutually exclusive while  $j=1,2,3$ and  $u=0,1$. Let
$$ 
\beta_{j+3u,\varphi_{j,u}}\ne \pm1, \, j=2,3.$$
For  $\theta_{1,u}>\varphi_{1,u}+3,
	\quad \theta_{2,u}=\varphi_{2,u}+3,
	\quad \theta_{3,u}=\varphi_{3,u}+3,$ we have:

  	\bigskip
  	
  	\noindent 1)  For  Subcase$^{(u)}$ 8.1: $\theta_1=2$, $\varphi_2=0$, $\varphi_3=1$,
  	 
  	 \smallskip
  
	\noindent  $\det\Gamma_{3u,k+3u}\neq 0$ if and only $\beta_{2+3u,\varphi_{2,u}}$ and $\beta_{3+3u,\varphi_{3,u}}$ satisfy the conditions of Lemma \ref{lem8.1}.

\smallskip

\noindent 2) For Subcase$^{(u)}$ 8.2:  $\theta_1=5$, $\varphi_2=0$, $\varphi_3=1$, 
	 
  	 \smallskip
  
$\det\Gamma_{3u,k+3u}\neq 0$ if and only $\beta_{2+3u,\varphi_{2,u}}$ and $\beta_{3+3u,\varphi_{3,u}}$ satisfy the conditions of Lemma \ref{lem8.2}.

\smallskip

\noindent 3) For Subcase$^{(u)}$  8.3: $\theta_1=1$, $\varphi_2=0$, $\varphi_3=2$,

	 \smallskip
  
$\det\Gamma_{3u,k+3u}\neq 0$ if and only $\beta_{2+3u,\varphi_{2,u}}$ and $\beta_{3+3u,\varphi_{3,u}}$ satisfy the conditions of Lemma \ref{lem8.3}.

\smallskip
\noindent 4) For Subcase$^{(u)}$  8.4: $\theta_1=4$, $\varphi_2=0$, $\varphi_3=2$,
	 \smallskip
  
$\det\Gamma_{3u,k+3u}\neq 0$ if and only $\beta_{2+3u,\varphi_{2,u}}$ and $\beta_{3+3u,\varphi_{3,u}}$ satisfy the conditions of Lemma \ref{lem8.4}.

\smallskip
\noindent 5) For Subcase$^{(u)}$  8.5: $\theta_1=0$, $\varphi_2=1$, $\varphi_3=2$,
 \smallskip
  
$\det\Gamma_{3u,k+3u}\neq 0$ if and only $\beta_{2+3u,\varphi_{2,u}}$ and $\beta_{3+3u,\varphi_{3,u}}$ satisfy the conditions of Lemma \ref{lem8.5}.

\smallskip
 \noindent 6) Subcase$^{(u)}$  8.6: $\theta_1=3$, $\varphi_2=1$, $\varphi_3=2$,
 \smallskip
  
$\det\Gamma_{3u,k+3u}\neq 0$ if and only $\beta_{2+3u,\varphi_{2,u}}$ and $\beta_{3+3u,\varphi_{3,u}}$ satisfy the conditions of Lemma \ref{lem8.6}.

\smallskip
  \end{proposition}
  \begin{remark}Since the $\varphi_j$, $j=1,2,3$ are mutually distinct, then   Case$^{(u)}$ 9: $\theta_{1,u}<\varphi_{1,u}+3,
	\quad \theta_{2,u}=\varphi_{2,u}+3,
	\quad \theta_{3,u}=\varphi_{3,u}+3$ is the particular Subcase$^{(u)}$  8.5  of Case$^{(u)}$  8, where $\theta_1=0$ is replaced by $\varphi_1=0$.  \end{remark}

 \noindent   Case$^{(u)}$ 10: $\theta_{1,u}=\varphi_{1,u}+3,
	\quad \theta_{2,u}=\varphi_{2,u}+3,
	\quad \theta_{3,u}=\varphi_{3,u}+3$,\\
	Since  $\theta_{j,u}=\varphi_{j,u}+3$, $j=1,2,3$ and $u=0,1$, then it follows from Assumption \ref{3.1.1} that $\varphi_{1,u}=0$, $\varphi_{2,u}=1$ and $\varphi_{3,u}=2$. Hence, \begin{equation}\label{det10}
	\gamma_{3u,k+3u}=
\begin{pmatrix}
	\phi_{1,\tau_{k,u}} &\phi_{1,\tau_{k,u}+2} &\phi_{1,\tau_{k,u}}\\
	\phi_{2,\tau_{k,u}} &\xi^{2}\phi_{2,\tau_{k,u}+2} &\xi^{4}\phi_{2,\tau_{k,u}}\\
	\phi_{3,\tau_{k,u}}
		&\xi^{4}\phi_{3,\tau_{k,u}+2} &\xi^{8}\phi_{3,\tau_{k,u}}
	\end{pmatrix}.\end{equation}
\begin{align}\label{Case 9.tmz1}
	\det\gamma_{3u,k+3u}&=(\xi^8-1)\phi_{3, \tau_{k,u}}(\xi^2\phi_{1, \tau_{k,u}}\phi_{2, \tau_{k,u}+2}-\phi_{2, \tau_{k,u}}\phi_{1, \tau_{k,u}+2})\nonumber\\&\qquad -(\xi^4-1)\phi_{2, \tau_{k,u}}(\xi^4\phi_{1, \tau_{k,u}}\phi_{3, \tau_{k,u}+2}-\phi_{3, \tau_{k,u}}\phi_{1, \tau_{k,u}+2})\nonumber\\&=\beta_{3+3u, \varphi_{3,u}}\Big[\beta_{2+3u, \varphi_{2,u}}(-i\beta_{1+3u, \varphi_{1,u}}(2\xi^4+1)+i^{\tau_{k,u}}(2\xi^2-1)) \nonumber\\&\quad+(-1)^{\tau_{k,u}}i(\xi^4+\xi^2)-3i^{\tau_{k,u}}\beta_{1+3u, \varphi_{1,u}}\Big]\nonumber\\&\quad+\sqrt3i^{\tau_{k,u}+1}(\beta_{1+3u, \varphi_{1,u}}\beta_{2+3u, \varphi_{2,u}}+(-1)^{\tau_{k,u}+1})\nonumber\\&\quad+(-1)^{\tau_{k,u}+1}i(3\sqrt3i\beta_{2+3u, \varphi_{2,u}}+(\xi^2+1)\beta_{1+3u, \varphi_{1,u}}-3\xi^2)
	\nonumber\\&=\beta_{3+3u, \varphi_{3,u}}\Big[\beta_{2+3u, \varphi_{2,u}}(\sqrt3\beta_{1+3u, \varphi_{1,u}}+\sqrt3i^{\tau_{k,u}+1}) \nonumber\\&\quad+(-1)^{\tau_{k,u}}-3i^{\tau_{k,u}}\beta_{1+3u, \varphi_{1,u}}\Big]\nonumber\\&\quad+\sqrt3i^{\tau_{k,u}+1}(\beta_{1+3u, \varphi_{1,u}}\beta_{2+3u, \varphi_{2,u}}+(-1)^{\tau_{k,u}+1})\nonumber\\&\quad+(-1)^{\tau_{k,u}+1}i(3\sqrt3i\beta_{2+3u, \varphi_{2,u}}+(\xi^2+1)\beta_{1+3u, \varphi_{1,u}}-3\xi^2).
\end{align}
Let \begin{align}\label{alphaCase10}
\alpha_{3+3u, \varphi_{3,u}}=&
(-1)^{\tau_{k,u}}-3i^{\tau_{k,u}}\beta_{1+3u, \varphi_{1,u}}+\sqrt3i^{\tau_{k,u}+1}\big(\beta_{1+3u, \varphi_{1,u}}\beta_{2+3u, \varphi_{2,u}}\nonumber\\&\quad+(-1)^{\tau_{k,u}+1}\big)\nonumber\\&\quad+(-1)^{\tau_{k,u}+1}i(3\sqrt3i\beta_{2+3u, \varphi_{2,u}}+(\xi^2+1)\beta_{1+3u, \varphi_{1,u}}-3\xi^2)
\end{align}
and 
\begin{align}\label{deltaCase10}
\delta_{3+3u, \varphi_{3,u}}=\beta_{2+3u, \varphi_{2,u}}(\sqrt3\beta_{1+3u, \varphi_{1,u}}+\sqrt3i^{\tau_{k,u}+1}) +(-1)^{\tau_{k,u}}-3i^{\tau_{k,u}}\beta_{1+3u, \varphi_{1,u}}.
\end{align}
Observe that $\det\gamma_{3u,k+3u}\ne 0$ if $\beta_{2+3u, \varphi_{2,u}}=\dfrac{3i^{\tau_{k,u}+1}\beta_{1+3u, \varphi_{1,u}}+(-1)^{\tau_{k,u}+1}}{\sqrt3(\beta_{1+3u, \varphi_{1,u}}+i^{\tau_{k,u}+1})}$. However, if  $\beta_{2+3u, \varphi_{2,u}}\ne\dfrac{3i^{\tau_{k,u}+1}\beta_{1+3u, \varphi_{1,u}}+(-1)^{\tau_{k,u}+1}}{\sqrt3\big(\beta_{1+3u, \varphi_{1,u}}+i^{\tau_{k,u}+1}\big)}$ then
 \begin{align*}\det\gamma_{3u,k+3u}\ne 0\Longleftrightarrow \beta_{3+3u, \varphi_{3,u}}\ne -\dfrac{\alpha_{3+3u, \varphi_{3,u}}}{\delta_{3+3u, \varphi_{3,u}}}.\end{align*}

\begin{proposition}\label{prop10}
	Let $\xi= \sqrt{3}/{2}+i/{2} $, $\theta_{j,u}\in\{-\infty,0,1,2,3,4,5\}$ and $\varphi_{j,u}\in \{-\infty, 0,1,2\}$,  where $\theta_{j,u}=p_{j+3u}$ and $\varphi_{j,u}=q_{j+3u}$ are mutually exclusive while  $j=1,2,3$ and  $u=0,1$. Let $\alpha_{3+3u, \varphi_{3,u}}$ and $\delta_{3+3u, \varphi_{3,u}}$  as defined in \eqref{alphaCase10} and \eqref{deltaCase10}. Let 
$$ 
\beta_{j+3u,\varphi_{j,u}}\ne \pm1, \, j=1,2,3.$$
For  $\theta_{1,u}=\varphi_{1,u}+3,
	\quad \theta_{2,u}=\varphi_{2,u}+3,
	\quad \theta_{3,u}=\varphi_{3,u}+3,$:
	$$\det\Gamma_{3u,k+3u}\neq 0 $$
if and only one of the following conditions holds:\\\\
1) $\beta_{2+3u, \varphi_{2,u}}=\dfrac{3i^{\tau_{k,u}+1}\beta_{1+3u, \varphi_{1,u}}+(-1)^{\tau_{k,u}+1}}{\sqrt3(\beta_{1+3u, \varphi_{1,u}}+i^{\tau_{k,u}+1})}$,\\\\
2) $\beta_{2+3u, \varphi_{2,u}}\ne\dfrac{3i^{\tau_{k,u}+1}\beta_{1+3u, \varphi_{1,u}}+(-1)^{\tau_{k,u}+1}}{\sqrt3(\beta_{1+3u, \varphi_{1,u}}+i^{\tau_{k,u}+1})}$ and $\beta_{3+3u,\varphi_{3,u}}\ne -\dfrac{\alpha_{3+3u, \varphi_{3,u}}}{\delta_{3+3u, \varphi_{3,u}}}$.
\end{proposition}

\noindent Thus, it follows from Proposition \ref{4.5.1}, Proposition \ref{prop5}, Proposition \ref{prop6}, Proposition \ref{prop7}, Proposition \ref{prop8}, Proposition \ref{prop10},   and \cite[Definition 7.3.1 and Proposition 4.1.7]{mennicken2003non} that:

\begin{theorem}
	The problems \eqref{3.1}, \eqref{3.2} are Birkhoff regular if and only if there are $r_{0}, r_{1}\in\{1,2,3,4,5,6,7,8,9\}$ such that the conditions $C(r_{0},u)$ and $C(r_{1},u)$ hold.
\end{theorem}

\bibliographystyle{amsplain}
\bibliography{birkhoff}

\end{document}